\documentclass[12pt]{article}
\usepackage{e-jc}

\usepackage{amsthm,amsmath,amssymb}

\usepackage{graphicx}

\usepackage[colorlinks=true,citecolor=black,linkcolor=black,urlcolor=blue]{hyperref}

\theoremstyle{plain}
\newtheorem{theorem}{Theorem}

\newtheorem{corollary}[theorem]{Corollary}
\newtheorem{proposition}[theorem]{Proposition}

\newtheorem{observation}[theorem]{Observation}

\theoremstyle{definition}
\newtheorem{definition}[theorem]{Definition}
\newtheorem{example}[theorem]{Example}

\theoremstyle{remark}
\newtheorem{remark}[theorem]{Remark}

\def\BC{\mathbb C}

\def\E{\mathcal E}
\def\D{\Delta}

\def\G{\mathcal G}

\def\F{\mathcal F}

\def\P{\mathcal P}

\def\E{\mathcal E}

\def\a{\alpha}

\def\S{\Sigma}

\def\g{\gamma}
\def\k{\kappa}

\def\sign{\mathrm{sign}}

\def\Pfaf{\operatorname{Pfaf}}
\def\Arf{\mathrm{Arf}}

\title{\bf Discrete Dirac Operators, Critical\\ Embeddings and Ihara-Selberg Functions}

\author{Martin Loebl\thanks{Partially supported by the Czech Science Foundation under the 
contract number P202-13-21988S.}\\ 
\small Dep.~of Applied Mathematics\\[-0.8ex]
\small Malostranske n. $25$, 118 00 Praha 1 \\[-0.8ex] 
\small Charles University, Czech republic\\
\small\tt loebl@kam.mff.cuni.cz\\
\and
Petr Somberg\thanks{Supported by the grant GACR P201/12/G028.}  \\
\small Mathematical Institute\\[-0.8ex]
\small Sokolovska 83, 180 00 Praha 8 \\[-0.8ex]
\small Charles University, Czech republic\\
\small\tt somberg@karlin.mff.cuni.cz 
}

\date{\dateline{Sep 24, 2013}{Dec 18, 2014}\\
\small Mathematics Subject Classifications: 05C50, 82B20, 57M15}

\begin{document}

\maketitle

\begin{abstract}
  The aim of the paper is to formulate a discrete analogue of the claim made by 
Alvarez-Gaume et al., (\cite{a}), realizing the partition function of the 
free fermion
on a closed Riemann surface of genus $g$ as a linear combination of 
$2^{2g}$ Pfaffians of Dirac operators.
Let $G=(V,E)$ be a finite graph embedded in a closed Riemann surface $X$ of genus $g$, 
$x_e$ the collection of independent variables associated with each edge $e$ of $G$ 
(collected in one vector variable $x$) and $\S$ the set of all 
$2^{2g}$ Spin-structures on $X$.  We introduce $2^{2g}$ rotations $rot_s$ and  
$(2|E|\times 2|E|)$ matrices $\D(s)(x)$, $s\in \S$, of the transitions between 
the oriented edges of $G$ determined by rotations $rot_s$.
 We show that the generating function for the even subsets of edges of $G$, i.e., 
the Ising partition function, is a linear combination of the square roots of $2^{2g}$ 
Ihara-Selberg functions $I(\D(s)(x))$ also called Feynman functions. 
By a result of Foata--Zeilberger holds $I(\D(s)(x))=\det(I-\D'(s)(x))$, where 
$\D'(s)(x)$ is obtained from $\D(s)(x)$ by replacing some entries by $0$.
Thus each Feynman function is computable in polynomial time. 
We suggest that in the case of critical embedding of a bipartite graph $G$,
the Feynman functions provide suitable discrete analogues for the Pfaffians 
of discrete Dirac operators.

  \bigskip\noindent \textbf{Keywords:} discrete conformal structure, 
	critical embedding, Ihara-Selberg function, discrete Dirac operator, 
	Ising partition function
\end{abstract}


\section{Introduction}
\label{sec.int}
It is well known, cf. \cite{m}, \cite{k} and references therein, how to formulate 
the notion of a critical embedding 
of a finite graph in a closed Riemann surface 
in such a way that one can read the critical values 
of independent variables (i.e., the edge weights or coupling constants in the physical terminology) 
attached to the edges of the Dimer and the Ising problems on $G$ out of the collection of angles of 
the embedding termed discrete conformal structure. 
It is rather attractive task to study whether some of the properties associated to the 
notion of criticality in statistical physics may also be derived from the collection of 
geometric data attached to critical embeddings.

The main theme of the present paper is the formulation of certain discrete analogue of the  
claim made by Alvarez-Gaume et all, see \cite{a}, that the partition function of free fermion on a closed Riemann surface 
of genus $g$ is a linear combination of $2^{2g}$ Pfaffians of Dirac operators.
The theory of free fermion is generally accepted to be closely related to the criticality of both the Dimer and the Ising 
problems on graphs. 

The Dimer problem and its determinant type solution are of considerable continual
interest, see e.g., \cite{k}, \cite{c}. In \cite{bt}, the authors study the critical Ising model by reducing it 
via the determinant type method to the Dimer model. 
The case of planar graphs is well-understood in this setting, \cite{k}.
However, as proved in \cite{c} and explained in detail in the Appendix, Subsection \ref{sub.kc},
if one wants to obtain a discrete analogue of the claim in \cite{a} for the Dimer
model in a surface of a positive genus, one has to take into account global restrictions 
on the graph embedded into the Riemann surface. In particular, the conditions on Kasteleyn flatness 
in Corollary \ref{cor.o} need to be satisfied.
For the sake of completeness, the combinatorial approach 
used to describe the determinant type reduction of the critical Dimer model on Riemann surfaces of 
positive genus is briefly overviewed in the Appendix.

The present paper investigates the question whether, using the geometric data provided by a critical embedding of a graph $G$
in a Riemann surface $X$ and the weights of the edges of $G$ in such a way that  
the matrix of these weights captures the basic properties of the discrete Dirac operator, the partition function of the Dimer and 
Ising models may be evaluated explicitly.

In particular, we propose to overcome the previously mentioned limitations and obstacles for the determinant-type 
reasoning, and replace the determinants by so called Ihara-Selberg functions of the graph $G$. This in turn allows
to rewrite the generating function of the even subsets of edges of $G$ (the Ising model partition function associated to the graph $G$) as a linear combination of $2^{2g}$ square roots of Ihara-Selberg functions which we call Feynman functions. We build the Feynman functions in a combinatorial way in order to capture the emergence of the rotations $rot_s$ for Spin structures $s$ in the analysis of the Ising partition function, where the realization of $rot_s$ exploits the space of quadratic forms on $H_1(X,F_2)$. This approach makes use of the 
original treatment of the Ising partition function via Ihara-Selberg functions by Sherman, see \cite{s}, \cite{l}. 
Our basic notion of a $g-$graph $G^g$ allows rather convenient treatment of
a planar model for the embedding of the graph $G$ in a closed Riemann surface of genus $g$. 
The $g-$graphs have been successfully used recently in \cite{lm} in a related situation. 
Moreover, the planar model may improve the insight into the analysis of the limiting processes, which
has not been done yet. The rotations are treated in a different way by Johnson \cite{j}.


\subsection{Critical embeddings of graphs into closed Riemann surfaces}
\label{sub.critical}
In the present subsection we review the notion of critical embedding associated to 
discrete conformal structure on a graph embedded in a Riemann surface, \cite{m}.

Let us consider a triple $(X,G,\varphi)$, where $X$ is a closed Riemann surface, $G=(V,E)$ 
a graph and 
$\varphi : G\hookrightarrow X$ an embedding. The image of $\varphi$ defines a CW decomposition
of $X$ such that $X\setminus \varphi(G)$ is a disjoint union of open faces.
A useful local description of $X$ is given by local charts 
$\{\varphi_j:U_j\to\BC|\bigcup _jU_j=X\}_j$ covering $X$, equipped with the flat metric and
finite number of conical singularities at $\{P_j\}_j$  constrained by Gauss-Bonnet formula.

\begin{figure}[h]
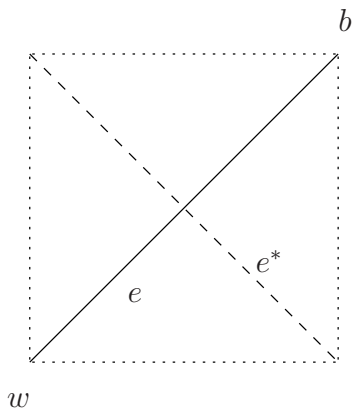

\begin{center}
\input diamond.pstex_t
\end{center}
\caption{Diamond $F(e)$ of edge $e$.}
\label{fig.diamond}
\end{figure}

An embedding of a graph $G$ in a surface $X$ induces the dual graph $G^*$ embedding. 
We regard
$G^*$ as an abstract graph with natural embedding into $X$, so that each vertex of $G^*$ lies
on the face of the embedding of $G$ it represents. The central notion related to the couple $G,G^*$ 
is that of the diamond graph - for simultaneous embedding of $G$ and $G^*$, 
the diamond graph $G^+$ has the vertex set equal to $V(G)\cup V(G^*)$ and the edges connecting the 
end-vertices of each dual pair of edges $e,e^*$ into a facial cycle $F(e)$ of $G^+$, which is a 
$4-$gon called a diamond (see Figure \ref{fig.diamond}).

An embedding $\varphi$ of $G^+$ is {critical} if each of its faces $F(e), e\in E$, 
is a { rhombus}, i.e., the following  conditions hold true with respect 
to the induced conformal class of metrics in a given local chart on $\varphi(X)$:
\begin{enumerate}
\item
The diagonals of each rhombus are perpendicular,
\item
The lengths of sides of all rhombi are the same.
\end{enumerate}
Notice that the first property is independent of the choice 
of local chart, because the transition maps $\varphi_i\circ\varphi_j^{-1}$ 
are conformal and so angle preserving, while the second condition 
depends on the choice of a representative metric in the conformal 
class.

\subsection{Dimer and Ising partition functions}
\label{sub.frfr}

In the present subsection we discuss the relationship between the Dimer and
the Ising partition functions with emphasis on their combinatorial description.

 Let us associate an independent variable $x_e$ with each edge $e\in E$ of $G$.  A subset $E'\subset E$ of edges is called {perfect matching} or {dimer arrangement}, if the induced graph $(V,E')$ has each vertex of degree one. Let $\P(G)$ denote the set of all perfect matchings of $G$.
We define the dimer partition function of $G$ by
\begin{eqnarray}
\P(G,x)= \sum_{M\in \P(G)}\prod_{e\in M}x_e,
\end{eqnarray}
where $x=(x_e)_{e\in E}$ is the vector of edge weights.

A subset $E'\subset E$ of edges is called { even} if the induced graph $(V,E')$ has each vertex of even degree. We denote by $\E(G)$ the set of even subsets of edges of $G$.

The generating function of the even sets of edges of $G$ is defined by
\begin{eqnarray}
\E(G,x)= \sum_{E' \in \E(G)}\prod_{e\in E'}x_e.
\end{eqnarray}
It is well known (see, e.g., \cite{lm}) that $\E(G,x)$ is equivalent to the Ising partition function on $G$ defined by
\begin{eqnarray}
Z_G^{\mathrm{Ising}}(\beta):= Z^{\mathrm{Ising}}_G(x)\Big|
  _{\textstyle{x_e:=e^{\beta J_e} \ \forall  e\in E}}\quad ,
\end{eqnarray}
 where $J_e$ $( e\in E)$ are the weights (coupling constants)
 associated with edges of the graph $G$, $\beta$ the 
 scale (inverse temperature) and
\begin{eqnarray}
Z_G^{\mathrm{Ising}}(x)= \sum_{\sigma:V\rightarrow \{1,-1\}} \  \prod_{e= \{u,v\}\in E}x_e^{\sigma(u)\sigma(v)}.
\end{eqnarray}

\subsection {The aim and scope of the article and its main results}
\label{sub.mr}

We shall close the first section by a brief formulation and overview of the main result
in the article.  

We denote by $\S$ the set of Spin structures on a Riemann surface $X$.
Following \cite{l}, \cite{j}, \cite{lm}, we associate 
the rotation to each Spin structure.
This enables us to define $(2|E|\times 2|E|)$ matrices $\D(s)(x)$, $s\in \S$, in a way that
each $\D(s)(x)$ is the transition matrix between the oriented edges determined by the rotation 
$rot_s$ corresponding to the Spin structure $s$.
The main result of our article is

\begin{theorem}
\label{thm.main}
Let $G$ be a graph embedded in a closed Riemann surface $X$ of genus $g$ 
and let $\S$ be the 
set of Spin structures on $X$.  
Then $\E(G,x)$ is a linear combination of $2^{2g}$ Feynman functions 
$F(\D(s)(x))$ for $s\in \S$ and $\D(s)$ the transition matrix of $G$. 
This is an Arf invariant formula and each 
$F(\D(s)(x))$ is computable in polynomial time. 

In addition, let $G$ be a 
bipartite graph critically embedded in the Riemann surface $X$ of genus $g$. 
Then for each Spin structure $s\in \S$ the transition matrix $\D(s)$ 
may be constructed out of the matrix of the discrete Dirac operator associated to $s\in \S$. 
\end{theorem}

The structure of our paper is organized as follows. In Section $2$ we introduce 
our main tool, the discrete Ihara-Selberg functions. We further construct a 
useful combinatorial model of closed Riemann surface $X$ of genus $g$, which 
enables to define $4^g$ rotations of prime reduced cycles on a graph critically
embedded in $X$. In Section $3$ we introduce quadratic forms, relate them to 
rotations and define the Feynman functions. In the end of this section we prove 
the first part of Theorem \ref{thm.main}, namely the Arf-invariant formula.
Section $4$ treats the topic of discrete Dirac operators, and concludes the 
proof of Theorem \ref{thm.main}. In the Appendix we summarize, for the reader's
convenience, the Pfaffian method and its limitations.
 
\subsection*{Acknowledgements}
The authors would like to thank David Cimasoni for extensive discussions about 
discrete Dirac operators and Gregor Masbaum for extensive discussions 
and suggestions about rotations and Spin structures.


\section{Discrete Ihara-Selberg functions and rotations}
\label{s.ciz}
In this section we assume that the graph $G$ is embedded in a closed Riemann surface $X$ of genus $g$. We suggest 
to consider the Ising partition function $\E(G,x)$ instead of the Dimer partition function
$\P(G,x)$, and the square root of certain Ihara-Selberg functions on the graph $G$ (which we call the Feynman functions) 
instead of the determinant. Similar, but much less advanced structure can be found in e.g., \cite{sm}, under the notion of Ising preholomorphic observable.

\subsection{Discrete Ihara-Selberg functions}
\label{sub.is}
Let $G=(V,E)$ be a graph. For $e\in E$ we denote by $o_e$ an orientation of $e$,
and $o_e^{-1}$ the reversed directed edge to $o_e$. As above, let $x= (x_e)_{e\in E}$ 
be the formal variables associated with edges of $G$. If $o$ is any orientation of the edge $e$, 
we associate the new variable
$x_o$ with it and always let $x_o= x_e$.

We consider an equivalence relation on the set of finite-length sequences $(z_1, \ldots, z_n)$
satisfying $z_1= z_n$: any such sequence is equivalent with each of its cyclic shifts.
The equivalence classes will be called {\em circular sequences}.

A circular sequence $p=v_1,o_1,v_2,o_2,...,o_n,v_{n+1}$ with $v_{n+1}=v_1$ and $o_i=(v_i,v_{i+1})$ is called a
{\it prime reduced cycle} if the following
conditions are satisfied:
$o_i \in \{o_e, o_e^{-1}: e\in E\}$, $o_i\neq o_{i+1}^{-1}$ and $(o_1,...,o_n)\neq Z^m$
for some sequence $Z$ and a natural number
 $m>1$.  We say that the ordered pair $(o_i, o_{i+1})$ is a transition of $p$ at $v_{i+1}$,
and $(o_n, o_1)$ is a transition of $p$ at $v_1$. We denote by $p^{-1}$ the prime reduced cycle which is the inverse of $p$.

\begin{definition} 
\label{def.z}
Let $G=(V,E)$ be a graph and assume that the vertex set $E$ is linearly ordered.
Let $M$ be $2|E|\times 2|E|$ matrix with entries $m(o,o')$,  $o,o'\in \{o_e, o_e^{-1}: e\in E\}$, where we think of $m(o,o')$ as 
the weight of the transition between directed edges $o,o'$ of $G$. If $p$ is a prime reduced cycle then we set 
$M(p)= \prod_{(o,o')\text{ a transition of }p}m(o,o')$.

We denote the set of prime reduced cycles
of $G$ by $\G$.
 The Ihara-Selberg function\index{function!Ihara-Selberg} associated to $G$ is defined by
\begin{eqnarray}
I(M)= \prod_{\gamma\in \G}(1-M(\gamma)),
\end{eqnarray}
where the infinite product is determined by the formal power series
\begin{eqnarray}
\prod_{\gamma\in \G}(1-M(\gamma))= \sum_{\F}(-1)^{|\F|}\prod_{\gamma \in \F}M(\gamma)
\end{eqnarray}
and the sum is over all finite subsets $\F$ of $\G$.

\end{definition}

A theorem of Foata--Zeilberger (see \cite{fz}), generalizing the seminal  theorem of Bass (see \cite{b}), states:

\begin{theorem}
\label{thm.fz}
We have
\begin{eqnarray}
I(M)= \det (I- M'),
\end{eqnarray}
where $M'$ is the matrix obtained from $M$ by letting $m'(o,o')= 0$ if $o'= o^{-1}$ and $m'(o,o')= m(o,o')$ otherwise.
\end{theorem}

\subsection{Combinatorial model of closed Riemann surfaces of genus $g$}
\label{sub.t}
 In this section we  restrict ourselves to the following standard representation of closed Riemann surface $X$ of genus $g$: 
we regard $X$ as a regular $4g-$gon $R$ (called the { base polygon}) in the plane with sides denoted anti-clockwise  by $z_1,\ldots, z_{4g}$,  and  the pairs of sides $z_i, z_{i+2}^{-1}$ and $z_{i+1}, z_{i+3}^{-1}$, $i=1, 5, \ldots, 4(g-1)+1$, identified. This defines a flat metric on $X$ 
with one conical singularity of angle $2\pi(2g-1)$. We restrict to embeddings of graphs on $X$ which meet the boundary of $R$ transversely. 
Moreover, to simplify the arguments, we only consider 
embeddings of graphs where each edge is represented by a straight line on $X$. The general embeddings may be treated analogously.

We now describe how an embedding of a graph in a Riemann surface can be used
to make its planar drawing of a special kind, cf. \cite{l}, \cite{lm}. 

\begin{definition}
\label{def.highway}
The {\em highway surface} $S_g$ consists of the base polygon $R$ and the bridges $R_1, \ldots, R_{2g}$, where 
\begin{enumerate}
\item
Each odd bridge $R_{2i-1}$ is a 
rectangle
with vertices $x(i,1), \ldots, x(i,4)$ numbered anti-clockwise. The bridge is glued
to $R$ so that its edge $[x(i,1),x(i,2)]$ is identified with the edge $[z_{4(i-1)+1},z_{4(i-1)+2}]$
and the edge $[x(i,3),x(i,4)]$ is identified with the edge $[z_{4(i-1)+3},z_{4(i-1)+4}]$.
\item
Each even bridge $R_{2i}$ is a 
rectangle
with vertices $y(i,1), \ldots, y(i,4)$ numbered anti-clockwise. It is glued
with $R$ so that its edge $[y(i,1),y(i,2)]$ is identified with the edge $[z_{4(i-1)+2},z_{4(i-1)+3}]$
and the edge $[y(i,3),y(i,4)]$ is identified with the edge
$[z_{4(i-1)+4},z_{4(i-1)+5}]$ (the indexes are always considered modulo $4g$.)
\end{enumerate}
\end{definition}

There is an orientation-preserving 
immersion $\Phi$ of $S_g$
into the complex plane which is injective except that for each $i=1, \ldots
g$, the images of the bridges
$R_{2i}$ and $R_{2i-1}$ intersect in a square (see Figure \ref{fig.torus}).


Now assume the graph $G$ is {piece-wise linearly embedded} into a closed Riemann surface $X$
of genus $g$. We realize $X$ as the union of $S_g$ and an additional disk
$R_\infty$ glued to the boundary of $S_g$. By an isotopy of the embedding
we may assume that $G$ does not meet the disk $R_\infty$ and all
vertices of $G$ lie in the interior of $R$. We may also assume that
the intersection of $G$ with any of the rectangular bridges $R_i$
consists of disjoint straight
lines connecting the two sides 
of $R_i$ glued to the base polygon $R$. The composition of the embedding of $G$
into $S_g$ with the immersion $\Phi$ yields a drawing $\varphi$ of
$G$ in the plane, where each edge of $G$ is represented by a piece-wise linear curve 
(see Figure \ref{fig.torus}).  A planar drawing of $G$ obtained in this way will be
called a {$g-$graph} and denoted by $G^g$.
Observe that double points of a $g-$graph can only come from the intersection of the images of
bridges under the immersion $\Phi$ of $S_g$ into the plane. Thus every double
point of a special drawing lies in one of the squares
$\Phi(R_{2i})\cap \Phi(R_{2i-1})$.

\begin{definition}
\label{def.rrr}
Let $G$ be embedded in $S_g$ and let $e$ be an edge of $G$. By definition, the embedding of $e$ intersects each bridge $R_i$  in disjoint straight lines. The number of these lines is denoted by $r_i(e)$. For a set $A$ 
of edges of $G$ we denote $r(A)$ be the vector of length $2g$ defined by 
$r(A)_i= \sum_{e\in A}r_i(e)$.
\end{definition}

\begin{definition}
\label{def.nnn}
Let $p$ be a prime reduced cycle of $G$. Then $p^g$ denotes the image of $p$
in $G^g$. 
\end{definition}

Clearly, $p\mapsto p^g$ gives bijective correspondence between the prime reduced cycles of $G$ 
and the prime reduced cycles of $G^g$.

\begin{figure}[h]
\begin{center}
\input torus.pstex_t
\end{center}
\caption{Immersions of edges $\{a,b\}$, $\{c,d\}$ and $\{e,f\}$ which cross a side of $R$.}
\label{fig.torus}
\end{figure}


\subsection{Rotations and the Ising partition function}
\label{sub.ff}
Let $G$ be a graph embedded in a closed Riemann surface $X$ of genus $g$
and let $G^g$ be its $g-$graph. We recall that each edge of $G$ in $G^g$ is 
represented by a piece-wise linear curve. Let $p$ be a prime reduced cycle of $G$ and
let $p^g$ the corresponding prime reduced cycle of $G^g$. We shall introduce $4^g$ 
{rotations} corresponding to $p^g$.

\begin{figure}[h]
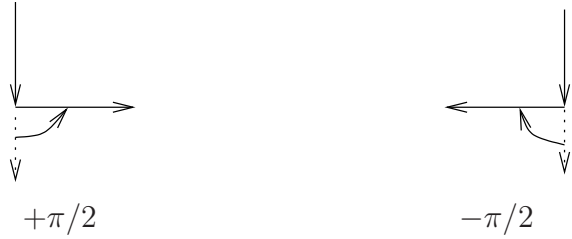

\begin{center}
\input tra.pstex_t
\end{center}
\caption{The angle of the transition.}
\label{fig.tra}
\end{figure}

First of all, we denote by $0$ the $0-$vector of length $2g$ and in analogy with 
the usual definition of the rotation of a regular closed curve in the plane we set
$rot_0(p^g)= \sum_t y_0(t)\pmod 2$, where we sum over all transitions of the linear
components of $p^g$.
If the transition $t$ consists in passing from directed segment $e$ to directed segment $e'$ then
$y_0(t)= z_0(t)(2\pi)^{-1}$, where $z_0(t)$ is the angle of the transition. The angle $z_0(t)$ is negative if the 
transition is clockwise, and $z_0(t)$ is positive if the transition is anti-clockwise (see Figure \ref{fig.tra}). 

Let $F_2$ be the field with two elements. Then we define $rot_s(p^g)$ for each arithmetic 
vector $s\in F_2^{2g}$ in the following way: if the transition $t$ consists in passing from the directed 
segment $e$ to the directed segment $e'$ and $e$ belongs to the immersion of bridge $R_i$ with 
$s_i=1$, we set $y_s(t)= y_0(t)+ 1$ and $y_s(t)= y_0(t)$ otherwise. 
Let $rot_s(p^g)= \sum_t y_s(t)$. We have $rot_s(p^g)= sr(p)+rot_0(p^g)\pmod 2$, where 
$sr$ denotes the scalar product of vectors $s$, $r$. Observe that for each $s$,  
$(-1)^{rot_s(p^g)}= (-1)^{rot_s((p^g)^{-1})}$. 

\begin{example}
Let $G$ be a toroidal graph embedded in the highway surface $S_1$, and let $p$ 
be a prime reduced cycle of $G$ intersecting each of the two bridges $R_i, i=1,2$ 
in just one segment.
The corresponding prime reduced cycle $p^1$ has exactly one self-intersection and $|rot_0(p^1)|= 0$.
Consequently, $rot_{00}(p^1)= rot_{11}(p^1)= 0$ and $rot_{01}(p^1)= rot_{10}(p^1)= 1$.
\end{example}


\begin{definition}
\label{def.new}
We introduce the equivalence relation on the set of prime reduced cycles of $G$: 
we say that $p_1$ is equivalent to $p_2$ if $p_1= (p_2)^{-1}$. The set of equivalence
classes for this equivalence relation is denoted by $[\G]$.
Analogously, we introduce the equivalence relation on the set of prime reduced cycles of $G^g$: 
 we say that $p^g_1$ is equivalent to $p^g_2$ if $p^g_1= (p^g_2)^{-1}$ and the set of equivalence
classes is denoted by $[\G]^g$.
\end{definition}

The previous considerations allow to introduce the functions
\begin{eqnarray}
F(G,x,s)=
\sum_{[\F]^g}(-1)^{|[\F]^g|}\prod_{[\gamma] \in [\F]^g}(-1)^{rot_s(\gamma)}\prod_{e\in \gamma}x_e,
\end{eqnarray}
where the sum is over all finite subsets $[\F]^g$ of $[G]^g$. 
In the proof of the first part of Theorem \ref{thm.main} in subsection \ref{sub.arrff} we observe 
that the functions $F(G,x,s)$ are the Feynman functions introduced in Definition \ref{def.ff}.

The following theorem appears in \cite{l} or the book \cite{l1}. We assign to $s\in F_2^{2g}$ the sign
according to $\sign(s)= (-1)^{\sum_{i=1,3,\ldots 2g-1} s_is_{i+1}}$. 

\begin{theorem} 
\label{thm.three}
Let $G$ be a graph with each vertex of degree equal to $2$ or $4$, embedded into 
a closed Riemann surface $X$ of genus $g$.
Then
\begin{eqnarray}
{\E}(G,x)= 2^{-g}\sum_{s\in F_2^{2g}} \sign(s)F(G,x,s).
\end{eqnarray}
\end{theorem}

We remark that Theorem \ref{thm.three} can be easily extended to general graphs and will be used in the proof
of Arf invariant formula (see Subsection \ref{sub.arrff}).

\begin{theorem}
\label{thm.four}
Let $G$ be a graph embedded into a closed Riemann surface $X$ of genus $g$.
Then
\begin{eqnarray}
{\E}(G,x)= 2^{-g}\sum_{s\in F_2^{2g}} \sign(s)F(G,x,s).
\end{eqnarray}
\end{theorem}

\begin{proof}
We will construct a graph $G'$ along with its embedding into $X$ so that the degree of each vertex of $G'$
is equal to $2$ or $4$, such that there are two subsets $Z, O$ of $E(G')$ and a bijection
$f: E(G')\setminus (Z\cup O) \rightarrow E(G)$ inducing
$$
{\E}(G,x)=
{\E}(G',z)|_{z_e:=x_{f(e)} \text{ if } e\notin  Z\cup O; z_e:=0 \text{ if } e\in Z; z_e:=1 \text{ if } e\in O} 
$$
and, for each $s\in F_2^{2g}$,
$$
F(G,x,s)=
F(G',z,s)|_{z_e:=x_{f(e)} \text{ if } e\notin  Z\cup O; z_e:=0 \text{ if } e\in Z; z_e:=1 \text{ if } e\in O}. 
$$
Theorem \ref{thm.four} then follows from Theorem \ref{thm.three}.  We construct graph $G'$
in two steps. Let $OD$ denote the set of the vertices of $G$ of an odd degree.
We start with $O= Z= \emptyset$. 

{\bf Step 1.}  If $OD\neq \emptyset$ then it is a standard observation of the graph theory that $G$ has a set of edge-disjoint
paths so that each vertex of $OD$ is an end-vertex of exactly one of the paths, and all the end-vertices of the paths are among the elements of $OD$ . In particular, $|OD|$ is even. 
Let $P$ denote the set of the edges of these paths. We construct graph $G_1$ from $G$
by adding, for each edge $\{u,v\}\in P$, a path of length $3$ with end-vertices $u,v$
(see Figure \ref{fig.zero}). We further let $Z$ be the set of all the edges of $G_1\setminus G$. Note that $G_1$ has all degrees even. 

\begin{figure}[h]
\begin{center}
\input zero.pstex_t
\end{center}
\caption{ }
\label{fig.zero}
\end{figure}

{\bf Step 2.} If $v$ is a vertex of $G_1$ of an even degree bigger than $4$ then we modify $G_1$ by splitting the degree of $v$ by introducing the new splitting edge; the operation
can be read off the Figure \ref{fig.one}. We repeat this step until the resulting graph $G'= (V,E')$ has all degrees equal to $2$ or $4$.  Finally we let $f$ be the tautological injection of $E$
into $E'$.

\begin{figure}[h]
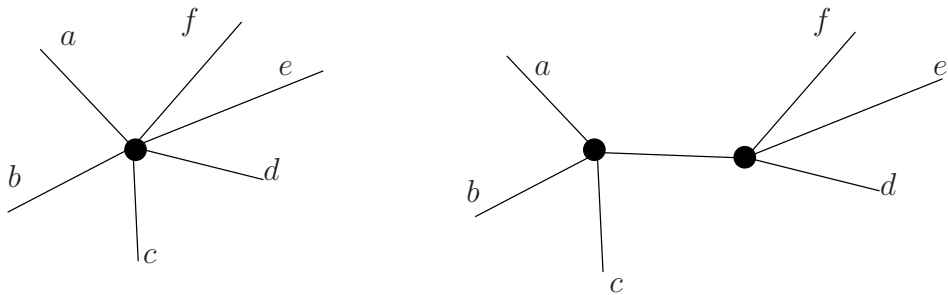

\begin{center}
\input one.pstex_t
\end{center}
\caption{ $O:=O\cup\{\text{ the new splitting edge }\}$}
\label{fig.one}
\end{figure}

It is straightforward to realize that all assumptions of the construction are satisfied after
application of finite number of these steps.

\end{proof}

In order to associate rotations to quadratic forms, we first need to study self-intersections of the prime reduced cycles which traverse each edge of the graph at most once.

\section{Feynman functions}

In the present section we relate the rotations of prime reduced cycles to
quadratic forms and moreover, we introduce the Feynman functions and prove 
a part of Theorem \ref{thm.main}: the Arf invariant formula.  

\subsection{Self-intersections of prime reduced cycles}
\label{s.ddpp}
\begin{figure}[h]
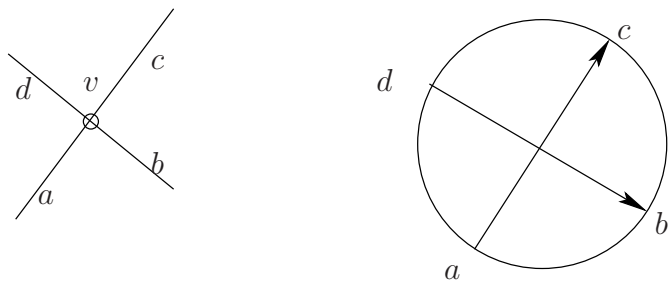

\begin{center}
\input crossing.pstex_t
\end{center}
\caption{Directed chord diagram $diag(p,v)$; $p= ac\ldots db\ldots a$}
\label{fig.crossing}
\end{figure}

Let $p$ be a prime reduced cycle of $G$ which traverses each edge of $G$ at most once. For each vertex $v$ of $G$, let $p(v)$ denote the set of the directed edges of $p$ incident with $v$, and  let $P(p,v)$ denote the partition of $p(v)$ into pairs which correspond to the {\em transitions} of $p$ at $v$. If a prime reduced cycle $p$ traverses each edge of $G$ at most once then the transitions of $p$ at $v$ are well described by the { directed chord diagram} $diag(p,v)$ (see Figure \ref{fig.crossing}):

\begin{definition}
\label{def.chord}
Let $p$ be a prime reduced cycle of $G$ which traverses each edge of $G$ at most once.
The directed chord diagram $diag(p,v)$ is obtained by taking the cyclic ordering of the edges of $G$ incident with $v$ and induced from the embedding of $G$ in $X$, and by introducing the directed chord $(e,e')$ for each class of $P(p,v)$ consisting of an orientation of $e$ followed by an orientation of $e'$.

We define the number of {\em self-intersections} of $p$ as the number of the pairs of intersecting chords of $diag(p,v)$, $v\in V$.
\end{definition}


We need to extend Definition \ref{def.chord} to the number of self-intersections of general prime reduced cycles, i.e., the prime reduced cycles that can go through an edge more than once. This can be done for instance in the following way. We define the infinite graph $\tilde{G}$ by replacing each edge $e$ by an infinite sequence of edges $e_1,\cdots, e_i, \cdots$
with the same end-vertex as $e$. We embed $\tilde{G}$ in $X$ so that we 'thicken' the embedding
of each edge $e$ of $G$, and embed the edges   $e_1,\cdots, e_i, \cdots$ to this thickened part of $e$ so that they are piece-wise
linear and internally disjoint. Next, for each prime reduced cycle $p$ of $G$ whose
circular sequence of directed edges is $(a_1\dots, a_k)$, $a_j$ being an orientation of edge $e(j)$ of $G$ (possibly $e(j)= e(l)$ for $j\neq l$),
we define prime reduced cycle $\tilde{p}$ in $\tilde{G}$ by replacing each $a_j$ by the same orientation of
$e(j)_j$.  It is important that the prime reduced cycle $\tilde{p}$ uses each edge of $\tilde{G}$
at most once. We thus define the number of the self-intersections of a prime reduced cycle $p$
of $G$ as the number of the self-intersections of the prime reduced cycle $\tilde{p}$ of $\tilde{G}$
(see Definition \ref{def.chord}). The number of the self-intersections of a prime reduced cycle $p^g$ of $G^g$ is defined analogously.

We note that for the prime reduced cycles of $G$ containing each edge of $G$ at most once this is consistent with Definition \ref{def.chord}. We also note that  the following basic property 
is satisfied.

\begin{observation}
\label{o.inf}
Let $p$ be a prime reduced cycle of $G$. Then we have $\pmod 2$
$$
1+ \text{ the number of the self-intersections of }p^g= rot_0(p^g).
$$
\end{observation}


\subsection{Rotations and Quadratic forms}
\label{sub.qf}

A subset of edges $E'\subset E$ of a graph $G$ embedded in a closed Riemann surface $X$ of genus $g$ is 
called even if each degree of the graph $(V,E')$ is even. The set of even subsets of $G$ is denoted by $\E(G)$.
Now we pass to the associated notion of {quadratic form}, cf. \cite{lm}.

Let $H:=H_1(X,F_2)$ be the first homology group of $X$
with coefficients in the field $F_2$. The construction of the highway surface gives us canonical 
basis of $H:=H_1(X,F_2)$, 
$$
a_1, b_1, \ldots, a_g, b_g,
$$
where $a_i$ corresponds to the class of the bridge $R_{2i-1}$ and $b_i$ corresponds to the class of $R_{2i}$,
$i=1,\ldots ,g$. Each element of $H$ is represented by the coordinate vector in $F_2^{2g}$, and two even 
subsets of edges 
$A, B$ belong to the same class in $H$ if and only if $r(A)= r(B)\pmod 2$ (see Definition \ref{def.rrr} for the definition of $r$). 

We recall that $H$ carries a non-degenerate skew-symmetric bilinear form
called intersection form and denoted by '$\cdot$'.
In the basis chosen above, it is determined by 
$a_i\cdot a_j= b_i\cdot b_j= 0$ and 
$a_i\cdot b_j= \delta_{ij}$    for all $i,j=1,\ldots, g$.

\begin{definition} 
\label{def.q}
 A {quadratic form on $(H,\cdot)$} associated to 
 the bilinear form '$\cdot$'
is a function $q:H\rightarrow F_2$ fulfilling
$q(x+y)=q(x)+q(y)+x\cdot y$ for all $ (x,y \in H)$.
\end{definition}
 We denote the set of quadratic forms on $H$ by $Q$. Each quadratic form is given by its values on 
the basis $a_1, b_1, \ldots, a_g, b_g$, and so the cardinality of $Q$ is $4^g$. 
 Let us denote by $q_0$ the quadratic form whose value on each of the basis vectors 
$a_i, b_i$ is zero. 
For each $z\in H$ let $q_z: H\rightarrow F_2$ be defined by $q_z(x)= q_0(x)+ zx$.
Then $q_z$ is a quadratic form and each quadratic form equals $q_z$ for some $z\in H$.

If $x\in H$ is a vector in $F_2^{2g}$ then we observe 
$\sign(x)= (-1)^{q_0(x)}$, where $\sign$ was introduced before the statement of Theorem \ref{thm.three}.

We now recall the definition of the {\em Arf invariant} of a quadratic form $q\in Q$. Let $N_0= 2^{g-1}(2^g+1)$ and $N_1= 2^{g-1}(2^g-1)$. Each quadratic form $q$ has either $N_0$ times value $0$ or $N_1$ times value $0$. In the first case we let $Arf(q)= 0$, and in the second case $Arf(q)= 1$.

The relevance of the Arf invariant in our considerations comes from the following fact, established as Lemma 2.10 in \cite{lm}: for each $x\in H$ 
\begin{align*}
2^{-g}\sum_{q\in Q}(-1)^{Arf(q)}(-1)^{q(x)}= 1.
\end{align*}

We will also use the following formula, established in the proof of Lemma 2.10., \cite{lm}: for 
each $z\in H$ holds
\begin{align}\label{signarf}
Arf(q_z)= q_0(z).
\end{align}

Notice that both quadratic forms $q_z$ and rotations $rot_z$ are parametrized by 
vectors $z\in F_2^{2g}$. Also recall that for a Riemann surface $X$ the Spin structures 
are given by the equivalence classes of square roots of the canonical bundle, and
these classes correspond to the $F_2$-valued first cohomology group $H^1(X,F_2)$. 
Because $X$ is assumed to be smooth, the Poincar\'e duality implies 
$H^1(X,F_2)\simeq H_1(X,F_2)$. 

There is a natural bijection between the set $\S$ of the 
Spin structures of $X$ and the set of quadratic forms on the $F_2-$valued 
first homology classes of $X$ (see \cite{j}). 
From now on we will consider the rotations $rot_s$ indexed by Spin structures on $X$. 

The usefulness of quadratic forms in studying rotations $rot_s$ was suggested to us by G. Masbaum, who also suggested Theorem \ref{thm.gr} below. 
\begin{theorem}
\label{thm.gr}
Let $p$ be a prime reduced cycle of $G$ and let $s\in \S$. Then 
\begin{eqnarray}
rot_s(p^g)= 1+ \text{ the number of the self-intersections of }p+ q_s(p)\, \pmod 2,
\end{eqnarray}
where $p^g$ is the realization of $p$ in the $g-$graph $G^g$.
\end{theorem}
\begin{proof}
By the definition of $rot_s$ it suffices to prove the statement for $s=0$. We have $\pmod 2$
$$
rot_0(p^g)= 1+ \text{ the number of the self-intersections of }p^g= 
$$
$$
1+ \text{ the number of the self-intersections of }p+ \sum_{i=1}^gr(p)_{2i-1}r(p)_{2i}=
$$
$$
1+ \text{ the number of the self-intersections of }p+ q_0(p),
$$
where the first equality follows by Observation \ref{o.inf} and the last one 
is a consequence of the definition $q_0$.

\end{proof}


\subsection{Feynman functions}
\label{sub.matt}
\begin{definition}
\label{def.dd}
Let $G$ be a graph embedded in a closed Riemann surface $X$ of genus $g$. Let $E^o= \cup\{o_e, o_e^{-1}: e\in E\}$, hence $|E^o|= 2|E|$. 
To each Spin structure $s\in \S$ we associate $\D(s)(x)$, the $|E^o|\times |E^o|$-matrix with entries $d_s(o,o')= (-1)^{y_s(o,o')} (-1)^{\k(o)}x_o'$, where 
\begin{enumerate}
\item
$y_s(o,o')$ were introduced in the beginning of the section \ref{sub.ff},
\item
$\k(o)= 0$ if $o$ is contained in the interior of $R$,
\item
$\k(o)= -3/4$ if the segment of the embedding of $o$ outside $R$ in $G^g$ is oriented oppositely to the anti-clockwise orientation of the boundary of $R$ (as directed edge $(a,b)$ in Figure \ref{fig.torus}), 
\item
$\k(o)= 3/4$ if the segment of the embedding of $o$ outside $R$ in $G^g$ is oriented in agreement with the anti-clockwise orientation of the boundary of $R$ (as directed edge $(a,b)$ in Figure \ref{fig.torus}). 
\end{enumerate}
We further define $\D'(s)(x)$ by declaring $d'_s(o,o')= 0$ if $o'= o^{-1}$ and $d'_s(o,o')= d_s(o,o')$ otherwise in $\D(s)(x)$.
\end{definition}
We are going to introduce the {\em Feynman functions}, see Definition \ref{def.ff} below.

\begin{definition}
\label{def.ff}
Let $s\in \S$ be a Spin structure on $X$. We define 
\begin{eqnarray}
F(\D(s)(x)):=
\sum_{[\F]}\prod_{[\gamma]\in [\F]}(-1)^{q_s(\gamma)+\text{ the number of self-intersections of } \gamma}\prod_{e\in \gamma}x_e.
\end{eqnarray}
We recall that $[\F]$ denotes the collection of all the finite sets of the equivalence classes of the reduced prime cycles of $G$. 
\end{definition}

\subsection{The Arf-invariant formula}
\label{sub.arrff}

The linear combination appearing in Theorem \ref{thm.main} is the Arf invariant formula,
\begin{eqnarray}\label{isarf}
{\E}(G,x)= 2^{-g}\sum_{s\in\S} (-1)^{\Arf(q_s)}F(\D(s)(x)).
\end{eqnarray}

\begin{proof} ({of the first part of Theorem \ref{thm.main}}, equation \ref{isarf}.)
 We have
\begin{eqnarray}
& & F(\D(s)(x))=
\sum_{[\F]}(-1)^{|[\F]|}\prod_{\gamma \in \F}(-1)^{1+q_s(\gamma)+\text{ the number of the self-intersections of } \gamma}\prod_{e\in \gamma}x_e
\nonumber \\ \nonumber
& & =\sum_{[\F]}(-1)^{|[\F]|}\prod_{\gamma \in \F}(-1)^{rot_s(\gamma^g)}\prod_{e\in \gamma}x_e = F(G,x,s),
\end{eqnarray}
where the second equality follows from Theorem \ref{thm.gr}.

Both $F(\D(s)(x))$ and $I(\D(s)(x))$ are written as an infinite product 
(see Definition \ref{def.z}). The factors of $I(\D(s)(x))$ are parametrized by 
the prime reduced cycles, while the factors of $F(\D(s)(x))$ are parametrized 
by the equivalence classes of prime reduced cycles. Hence each factor of $F(\D(s)(x))$
appears twice as factor of $I(\D(s)(x))$. 
Consequently, $F(\D(s)(x))$ is the square root of $I(\D(s)(x))$.
By Theorem \ref{thm.fz}, $I(\D(s)(x))= \det(I-\D'(s)(x))$ and each $F(\D(s)(x))$ is thus computable in polynomial time. 

Now the Arf-invariant formula follows from Theorem \ref{thm.four} and a simple observation done 
in equation \ref{signarf} mentioned in the end of subsection \ref{sub.qf}: for each $s$, $sign(s)= (-1)^{\Arf(q_s)}$.

\end{proof}

\section{Discrete Dirac operators}
\label{sub.DDo}
So far we considered arbitrary embeddings of a graph $G$ into a Riemann surface $X$. Now, let $G= (W,B,E)$ be a bipartite graph and 
the embedding critical. The vertices of $W$ are called {white} and the vertices of $B$ are called {black}.
We also assume in this subsection that the conical singularities are not located at the vertices of $G$. 

 We recall the notation $E^o= \cup\{o_e, o_e^{-1}: e\in E\}$.
Let $W^o$ be the subset of $E^o$ consisting of the edges directed from its black
vertices to its white vertices. Analogously, let $B^o= E^o\setminus W^o$ be the set of edges
directed from its white vertex to its black vertex. 
The key construction has the following structure: 

\begin{definition}
\label{def.t}
Let $T_G$ be the {\em directed transition graph}
of the orientations of edges of $G$, i.e. $V(T_G)= E^o$ and $(o,o')\in E(T_G)$ if
the head of $o$ is the tail of $o'$. 
\end{definition}

\begin{observation}
\label{o.kk1}
The graph  
$T_G$ is a directed bipartite graph, $T_G=(W^o,B^o, E(T_G))$, and the matrix  $\D(s)(x)$
is a weighted adjacency matrix of $T_G$.  
\end{observation}

The next observation follows directly from Definition \ref{def.dd}. 

\begin{observation}
\label{o.uhel}
Let $w^o\in W^o$ be an orientation of the edge $e$ entering the vertex $w\in W$. Let $b^o$ be a directed edge leaving $w$ and entering vertex $b\in B$.  Then  the entry $\D(s)(x)(w^o,b^o)$ equals 
$\g(s,w^o)e^{i\a(w^o,b^o)}x_{b_o}$, where 
$\a(w^o,b^o)= 1/2z_0(w^o,b^o)$ is half of the angle of the transition from $w^o$ to $b^o$; $\g(w^o)$ equals a complex number depending only on $s$ and $w^o$.
\end{observation}

 Let $l(e^*)$ denote the length of the dual edge $e^*$ of an edge $e$ of $G$.

\begin{definition}
\label{def.dada}
 We denote by $\D_2(s)$ the matrix obtained by taking the square of each entry of $\D(s)(x)$ and substituting the vector of the lengths of the dual edges for the squares of the variables: 
$(x_e)^2:= l(e^*), e\in E^o$.
\end{definition}

\begin{corollary} (of Observation \ref{o.uhel})
\label{c.dva}
If $w^0, w^1$ are two directed edges of $E^o$ entering the same vertex then the row of   $\D_2(s)$
indexed by $w^0$ is a complex multiple of the row of  $\D_2(s)$ indexed by $w^1$. 
\end{corollary}

The discrete Dirac operator $\mathcal{D}(s)$ of a bipartite graph $G= (V,E)$, $V= W\cup B$, corresponding to a Spin structure $s\in \S$, is defined as the weighted adjacency $(|V|\times|V|)-$matrix; \cite{k} contains the definition for the planar graphs and \cite{c} defines half of (or, the chiral part of) $\mathcal{D}(s)$, which uniquely determines $\mathcal{D}(s)$.
We observe in this subsection that the matrix $\D_2(s)$ is closely related to $\mathcal{D}(s)$. 

To that aim we fix the coordinate chart on the Riemann surface such that 
a black vertex $v\in B$ lies at $0\in\BC$ in this local chart. The directed 
edges emanating from a white vertex and terminating at $v$ are denoted
$e^v_1,\dots ,e^v_k$. Similarly, we fix an edge $e^v_w$ emanating from $v$ and
terminating at some white vertex $w\in W$. The $(vw)$-matrix element of the Dirac operator 
is $D_{(vw)}=l((e^v_w)^*)e^{i\alpha}$, where $\alpha$ is the angle between the edge $e^v_w$
and the real axis of the local coordinate chart. Note that 
the change in local coordinate chart does not preserve the real axis, but
does preserve the angle $\alpha$. On the other hand the $(e^v_i, e_w^v)$-entry of  
our matrix $\D_2(s)$ is equal to $l((e^v_w)^*)e^{i\alpha_i}$, where $\alpha_i$ is
the angle between $e^v_i,e_w^v$ measured anticlockwise, see Figure \ref{fig.tra}. 
In particular, in the local coordinate chart in which the 
edge $e^v_i$ lies on the real axis, we have $\alpha= \alpha_i$. 
Equivalently, the matrix coefficients in a 
given row of our matrix $\D_2(s)$ indexed by directed edge $e$ are (constant) multiples of the row in the Dirac matrix
corresponding to the terminal vertex of $e$. The same applies when 
one starts with the white vertex instead of the black one.

The basic observation about the matrices $\D_2(s)$ and $\mathcal{D}(s)$ can be formulated in the following way.

\begin{corollary}
\label{c.tri}
Let $w^o\in W^o$ be an orientation of the edge $e$ entering the vertex $w\in W$. Let $b^o$ be a directed edge leaving $w$ and entering vertex $b\in B$.  Then  $\D_2(s)(w^o,b^o)= c(w^o)\mathcal{D}(s)(w,b)$, where $c(w^o)$ is a complex number depending only on $w^o$. 
\end{corollary}

The last part of Theorem \ref{thm.main} is contained in the next Corollary.

\begin{corollary}
\label{c.ddo}
 Each matrix $\D_2(s)$, $s\in\Sigma$ a Spin structure, may be obtained from the  discrete Dirac operator 
$\mathcal{D}(s)$, whose chiral part is defined in \cite{c}, by finite number of operations:
 \begin{enumerate}
 \item
 Adding identical copy of a row, 
 \item
  Multiplying a row by a complex number,
\item
Adding $2|E|-|V|$ zero-entries to each row.
 \end{enumerate}
\end{corollary}

\begin{remark}
\label{r.conc}
The properties that graph $G$ is bipartite and critically embedded, and the condition that the conic singularities of the metric are not located at the vertices of $G$,  are needed in \cite{c} to show that constant functions belong to the kernel of the discrete Dirac operators. 
\end{remark}
Each Feynman function can be rewritten as the alternating sum of the traces of skew-symmetric powers of matrices  $\D'(s)(x)$,
which is an approach closely related to the origin of analytic torsion on Riemann surfaces. These considerations may lead to a realization 
of the analytic torsion as a continuous limit of the Feynman functions. 


\section{Conclusion}

The present paper is a contribution towards mathematical understanding of fermionic 
quantum field theory. In particular, we indirectly connect generating functions for 
the enumeration of even subsets 
of edges of finite graphs embedded in closed Riemann surfaces
with an expression for partition function of free fermion by means of Dirac operators.   
The crucial role is played by the discrete Ihara-Selberg functions. However, the question
of their limiting behavior remains elusive.


\section{Appendix: Pfaffian method, Dimers at criticality and Determinants}
\label{s.Dcd}

For the sake of completeness we briefly review the Pfaffian method and indicate its
global restrictions in order to obtain a discrete analogue of the result by Alvarez-Gaume
at all.

\subsection{Kasteleyn orientations}
\label{sub.Kor}
The first step in this theory is the reduction of $\E(G,x)$ to $\P(G_{\delta},x)$, where 
$G_{\delta}$ is obtained from $G$ by local changes which do not affect 
genus of the Riemann surface $X$ as a target for embedding of $G$. This construction 
(see e.g., \cite{lm}) is not relevant here, so we omit it and concentrate on the aspects of 
the Dimer partition function $\P(G,x)$.
 
Assume the vertices of $G$ are numbered
from $1$ to $n$. An orientation of $G$ is induced by prescribing one of the two 
possible directions to each edge of $G$.  If $D$ is an orientation of $G$, we denote 
by $A(G,D)$ the { skew-symmetric adjacency
matrix} of $D$ defined as follows: the diagonal entries of  $A(G,D)$ are zero,
and the off-diagonal entries are 
$$A(G,D)_{ij}=\sum_{e: i\to j} \pm x_e~,$$ where the sum is taken over all edges $e$
connecting vertices $i$ and $j$, and the sign in front of $x_e$ is $1$ if $e$ is
oriented from $i$ to $j$ in the orientation $D$, and $-1$
otherwise. It is well known that the Pfaffian
of $A(G,D)$ counts perfect matchings $\mathcal{P}(G)$ of the graph $G$ {with signs:} 
$$\Pfaf A(G,D) = \sum_{M\in \mathcal{P}(G)}\sign(M,D) \prod_{e\in
  M}x_e~,$$ where $\sign(M,D)=\pm 1$. We use this as the definition of the sign of a perfect
matching $M$ with respect to the orientation $D$.
The following statement is the basic result in the field:

\begin{theorem}[Kasteleyn \cite{ka}, Galluccio-Loebl \cite{gl},
  Tesler \cite{t}, Cimasoni-Reshetikhin \cite{cr}] \label{Pfaff-formula} If $G$
embeds into a Riemann surface of genus $g$, then there exist $4^g$
orientations $D_i$ ($i=1,\ldots, 4^g$) of $G$ such that the perfect matching
polynomial $\P(G,x)$ can be
expressed as a linear combination  of the Pfaffian polynomials
$\Pfaf A(G,D_i)(x)$.
Moreover, the orientations $D_i$ are { Kasteleyn orientations}, i.e.,
 they satisfy the following property: if $F$ is a facial cycle of $G$ then
 $F$ has an odd number of edges oriented in $D_i$ in agreement with the 
clockwise traversal of $F$.
\end{theorem}

The formula of Theorem \ref{Pfaff-formula} is called the {Arf-invariant formula}, as it is
based on the property of the Arf invariant for quadratic forms in
characteristic two. As far as we know, the relationship with the Arf
invariant was first observed in \cite{cr}. 
In \cite{ku}, Kuperberg introduced a generalization of the notion of Kasteleyn orientations, 
called {Kasteleyn flatness}.


\subsection{Kasteleyn curvature}
\label{sub.flat}
We recall that a graph $G=(V,E)$ is called { bipartite} if the set of the vertices $V$ may be partitioned 
into two sets $W, B$ so $|e\cap W|= |e\cap B|= 1$ for each $e\in E$.
In this and the next one subsection we restrict ourselves to finite bipartite graphs $G= (W,B,E)$, where $W,B$ 
are the two edge-less sets of vertices
of $G$ and $V= W\cup B$. We call the vertices in $W$ white and the vertices in $B$ black and
assume that $G$ has at least one perfect matching. In particular, we restrict to the case of 
equal cardinalities $|W|= |B|$.
By a cycle in $G$ we mean a subset of edges $C\subset E(G)$ which form a cycle. The cycle $C$ can be decorated
by one of the two possible orientations, i.e. one of the two possible ways of going around $C$. 
By an {oriented cycle} we mean a cycle decorated with an orientation. In this subsection we also use $O$ to denote the orientation of $G$ in which each edge is oriented from its black vertex to its white vertex.
 
\begin{definition}
\label{def.Kc}
Let $G$ be given with the weights $w(e), e\in E(G)$.
Let $C$ be an oriented cycle of $G$. We define 
$$
c(C)= (-1)^{|C|/2+1}\frac{\prod_{e\in C_+}w(e)}{\prod_{e\in C_-}w(e)}
$$
as the {\em Kasteleyn curvature} of $C$. Here $C_+$ denotes the subset of edges of $C$ whose orientation 
inherited from $C$ coincides with their orientation in $O$, and $C_-= C\setminus C_+$.
\end{definition}

\begin{definition}
\label{def.Kflat}
Let $G$ be a weighted graph, $w(e), e\in E(G)$, embedded in $X$. We say that $G$ is {\em Kasteleyn flat} 
if $c(F)=1$ for each face $F$
of the embedding, arbitrarily oriented. 
\end{definition}
Let $G=(W,B,E)$ be a bipartite graph equipped with the weights $w(e), e\in E(G)$. Let us fix a linear ordering on the set
$B\cup W$ such that the elements of $B$ precede the elements of $W$. 
This allows to introduce $D(w)$ as $(B\cup W)\times(B\cup W)$ skew-symmetric matrix defined 
by $D(w)_{uv}= w(uv)x_e$ if $u\in B$ and $e= uv\in E(G)$ resp. $D(w)_{uv}= -w(uv)x_e$ 
if $u\in W$ and $e= uv\in E(G)$, $D(w)_{uv}=0$ otherwise. We further
denote by $D_B(w)$ the $(B\times W)-$block of $D(w)$.
For $M$ a perfect matching of $G$ we denote by $t(M)$ the coefficient of $\prod_{e\in M}x_e$ in $\det(D_B(w))$.
Note that $t(M)= \sign(M)\prod_{e\in M}w(e)$, where $\sign(M)$ is the sign of $M$ relative to the 
fixed linear ordering on the set of vertices $B\cup W$.

The significance of flatness of the Kasteleyn curvature rests on the following observation (see \cite{ku}).

\begin{proposition}
\label{prop.two}
Let $C$ be a cycle of $G$ embedded in $X$ so that $C$ is a symmetric difference of two perfect 
matchings $M,N$ and the orientation of $C$ coincides with the orientation of the edges of $M$ in $O$.
Then 
$$
\frac{t(M)}{t(N)}= c(C).
$$
\end{proposition}

This result leads to 

\begin{proposition}
\label{prop.one}
If $G$ is embedded in the plane and its weight-function is Kasteleyn flat then $t(M)$ is constant.
\end{proposition}

We are able to prove a stronger statement. Let $G$ be embedded in a Riemann surface $X$ of genus $g$. We say
that the weight-function $w$ is {\em simple Kasteleyn flat} if it is Kasteleyn flat and moreover
$w(e)\in \{1,-1\}$ for each $e\in E$. The function $w$ may be viewed as assigning orientation $O(w)$ obtained
from orientation $O$ by reversing the orientation of all edges with negative weight. If $w$ is simple Kasteleyn flat 
the orientation
$D(w)$ is a { Kasteleyn orientation}, i.e. each face has an odd number of edges oriented clockwise.  
The next definition is motivated by \cite{cr}.

\begin{definition}
\label{def.eq}
We say that two weight-functions $w,w'$ are { equivalent} if $w$ can be obtained from $w'$ by finite number of { vertex multiplications},
where a vertex multiplication consists in choosing a vertex $v$ and a complex number $c\neq 0$ together with multiplication by $c$ the weights 
of all edges incident with $v$. 
\end{definition}

\begin{proposition}
\label{prop.dos}
Let $G$ be a graph embedded in a closed Riemann surface $X$ of genus $g$ and let $w$ be a Kasteleyn flat weight function,
satisfying in addition $c(C)\in \{1,-1\}$ for each cycle $C$. Then $w$ is equivalent to
a simple Kasteleyn flat weight-function $w'$, i.e., to a Kasteleyn orientation.
\end{proposition}
\begin{proof}
We may assume that $G$ is connected. Let $T=(V(G),E')$, $E'\subset E(G)$ be a spanning tree of $G$.
We can clearly perform vertex multiplications (in $G$) in a way the resulting weight function $w'$ 
satisfies $w'(e)= 1$ for each $e\in E'$.
Let $e\in E\setminus E'$. Then necessarily $w'(e)\in \{1,-1\}$ since $e$ forms a cycle (say, denoted by
$C$) with a subset of $E'$ and $c(C)$ has values in $\{1,-1\}$.
\end{proof}

Proposition \ref{prop.dos} along with Theorem \ref{Pfaff-formula} immediately imply

\begin{corollary}
\label{cor.o}
Let $G$ be a graph embedded in a closed Riemann surface $X$ of genus $g$ and $w$ be a Kasteleyn flat weight function
satisfying $c(C)\in \{1,-1\}$ for each cycle $C$. Then $\P(G,x)$ is a linear 
combination of $4^g$ determinants of  matrices obtained from $D_B(w)$ by multiplying some entries by $-1$.
\end{corollary}

\subsection{Kasteleyn curvature and criticality}
\label{sub.kc}
The computation of the Dimer partition function $\P(G,x)$ for a critically embedded bipartite graph $G$ is
based on the following strategy (see \cite{k}, \cite{c}): define the weight function $w$ 
using discrete geometric
information contained in the data of the critical embedding, and use this information to prove that the kernel 
of the corresponding matrix $D(w)$ 
has properties desired for a discrete analogue of the Dirac operator. Next, prove that $w$ is Kasteleyn flat and 
insert it into Corollary \ref{cor.o}. 

This approach works for planar bipartite graphs as well
(\cite{k}), since in this case each Kasteleyn flat weight function satisfies  $c(C)\in \{1,-1\}$ for each cycle $C$,
see \cite{ku}. 
On the other hand, the assumptions of Corollary \ref{cor.o} are quite restrictive 
for non-planar surfaces.
By Proposition \ref{prop.dos}, they are equivalent to $w(e)\in \{1,-1\}$ for each edge $e\in E$.
Moreover, it is shown in \cite{c} that the theory of Kasteleyn flatness can not go beyond
Corollary \ref{cor.o}, see genus $1$ example in \cite{c}, subsection 4.4. Examples.


\end{document}